\documentclass[leqno,12pt]{article}
\usepackage{amsmath,amssymb,rotating,a4wide,color}
\usepackage{wasysym}
\usepackage{hyperref}
\usepackage{adjustbox}

\usepackage[all,cmtip]{xy}

\usepackage{verbatim}

\newdir^{ (}{{}*!/-3pt/\dir^{(}}
\newdir_{ (}{{}*!/-3pt/\dir_{(}}

\entrymodifiers={+!!<0pt,\fontdimen22\textfont2>}

\def\Agemotext{\rotatebox[origin=c]{180}{$\textstyle\kern-1pt\varOmega$}}
\def\Agemoscript{\rotatebox[origin=c]{180}{$\scriptstyle\kern-1pt\varOmega\kern1pt$}}
\def\Agemoscriptscript{\rotatebox[origin=c]{180}{$\scriptscriptstyle\kern-1pt\varOmega\kern1pt$}}
\def\Agemo{{\mathchoice{\Agemotext}{\Agemotext}{\Agemoscript}{\Agemoscriptscript}}}

\newcommand{\rtuple}{{r}}

\def\strong{{\rm strong}}
\def\gr{\mathop{\rm gr}\nolimits}
\def\Sch{\mathop{\rm Sch}\nolimits}

\def\bigtimesdisplay{\mathop{\raise-2pt\hbox{\huge$\times$}}}
\def\bigtimestext{\mathop{\raise-1pt\hbox{\Large$\times$\kern-2pt}}}

\let\oldbigwedge\bigwedge
\def\newbigwedge{\mathord{\adjustbox{valign=B,totalheight=8.5pt}{$\oldbigwedge$}}}
\renewcommand{\bigwedge}{\newbigwedge}

\newbox\circbulletbox
\setbox\circbulletbox=\hbox{\,{\large$\circledcirc$}\kern-8.7pt\raise .6pt\hbox{$\bullet$}\,}

\let\le\leqslant
\let\ge\geqslant

\let\triangleleftnaked\triangleleft
\def\triangleleft{\mathrel{\triangleleftnaked}}
\def\depsilon{d^{\kern1pt\epsilon}}

\def\circVbig{\hbox{\text{\it\r{V}}}}
\def\circVscript{\hbox{\scriptsize\text{\it\r{V}}}}
\def\circVscriptscript{\mbox{\tiny\text{\it\r{V}}}}

\def\circVlimits_#1^#2{{\mathchoice%
   {\circVbig{}^{\kern2pt #2}_{\kern-2pt #1}}%
   {\circVbig{}^{\kern2pt #2}_{\kern-2pt #1}}%
   {\scriptstyle\circVscript{}^{\kern1.7pt #2}_{\kern-1pt #1}}%
   {\scriptscriptstyle\circVscriptscript{}^{\kern1.5pt #2}_{\kern-1pt #1}}%
   }}
\def\circVr_#1{\circVlimits_#1^r}
\def\circVs_#1{\circVlimits_#1^s}

\def\circWbig{\hbox{\text{\it\r{W}}}}
\def\circWscript{\hbox{\scriptsize\text{\it\r{W}}}}
\def\circWscriptscript{\mbox{\tiny\text{\it\r{W}}}}

\def\circWlimits_#1^#2{{\mathchoice%
   {\circWbig{}^{\kern2pt #2}_{\kern-2pt #1}}%
   {\circWbig{}^{\kern2pt #2}_{\kern-2pt #1}}%
   {\scriptstyle\circWscript{}^{\kern1.7pt #2}_{\kern-1pt #1}}%
   {\scriptscriptstyle\circWscriptscript{}^{\kern1.5pt #2}_{\kern-1pt #1}}%
   }}

\def\OM{\mathchoice
  {\rlap{\kern3.2pt$\overline{\phantom{L}}$}M}
  {\rlap{\kern3.2pt$\overline{\phantom{L}}$}M}
  {\rlap{\kern2.4pt$\scriptstyle\overline{\phantom{L}}$}M}
  {\rlap{\kern1.8pt$\scriptscriptstyle\overline{\phantom{L}}$}M}}

\def\mycirc{{\kern1pt\circ\kern2pt}}

\def\Aut{\mathop{\rm Aut}\nolimits}

\def\Gal{\mathop{\rm Gal}\nolimits}

\def\Mat{\mathop{\rm Mat}\nolimits}

\def\Lie{\mathop{\rm Lie}\nolimits}
\def\Ad{\mathop{\rm Ad}\nolimits}

\def\GL{\mathop{\rm GL}\nolimits}

\def\PGL{\mathop{\rm PGL}\nolimits}

\def\ab{{\rm ab}}

\let\phi\varphi

\def\theta{\vartheta}
\let\epsilon\varepsilon
\let\setminus\smallsetminus

\newtheorem{Thm}{Theorem}[section]
\newtheorem{Prop}[Thm]{Proposition}
\newtheorem{Lem}[Thm]{Lemma}

\newtheorem{Cor}[Thm]{Corollary}

\newtheorem{Rem}[Thm]{Remark}

\newtheorem{Ass}[Thm]{Assumption}

\numberwithin{Thm}{section}
\def\UseTheoremCounterForNextEquation{\setcounter{equation}{\value{Thm}}\addtocounter{Thm}{1}}

\def\qed{{\hskip0pt\unskip\unskip\nobreak\hfil\penalty50
          \hskip1em\hbox{}\nobreak\hfil
           {$\square$}
          \parfillskip=0pt\finalhyphendemerits=0
          \par}\medskip}

\newenvironment{Proof}
               {\noindent{\bf Proof.}\ }
               {\qed}

               {\noindent{\bf Proof of #1.}\ }
               {\qed}

\newcommand{\BF}{{\mathbb{F}}}

\newcommand{\BQ}{{\mathbb{Q}}}

\newcommand{\BZ}{{\mathbb{Z}}}

\newcommand{\CG}{{\cal G}}
\newcommand{\CH}{{\cal H}}

\newbox\mybox
\def\arrover#1{\mathrel{
       \setbox\mybox=\hbox spread 1.4em
              {\hfil$\scriptstyle#1$\hfil}
       \vbox{\offinterlineskip\copy\mybox
             \hbox to\wd\mybox{\rightarrowfill}}}}

\def\larrover#1{\mathrel{
       \setbox\mybox=\hbox spread 1.4em
              {\hfil$\scriptstyle#1\vphantom{g}$\hfil}
       \vbox{\offinterlineskip\copy\mybox
             \hbox to\wd\mybox{\leftarrowfill}}}}

\def\ontoover#1{\mathrel{
       \setbox\mybox=\hbox spread 1.4em
              {\hfil$\scriptstyle#1\vphantom{g}$\hfil}
       \vbox{\offinterlineskip\copy\mybox
             \hbox to\wd\mybox{\rightarrowfill\hskip-2.8mm
                               $\rightarrow$}}}}
\def\leftontoover#1{\mathrel{
       \setbox\mybox=\hbox spread 1.4em
              {\hfil$\scriptstyle#1\vphantom{g}$\hfil}
       \vbox{\offinterlineskip\copy\mybox
             \hbox to\wd\mybox{$\leftarrow$\hskip-2.8mm
                               \leftarrowfill}}}}
\let\longto\longrightarrow
\let\into\hookrightarrow
\let\onto\twoheadrightarrow

\def\isoto{\mathrel{
       \setbox\mybox=\hbox spread 0.9em
              {\hfil$\scriptstyle\sim$\hfil}
       \vbox{\offinterlineskip\copy\mybox
             \hbox to\wd\mybox{\rightarrowfill}}}}

\def\Bigskip{\bigskip\bigskip}

\begin{document}

\title{Schur $\sigma$-groups of type (3,3)}

\author{Richard Pink\\[12pt]
\small Department of Mathematics \\[-3pt]
\small ETH Z\"urich\\[-3pt]
\small 8092 Z\"urich\\[-3pt]
\small Switzerland \\[-3pt]
\small pink@math.ethz.ch\\[12pt]}

\date{May 14, 2025}

\maketitle

\Bigskip

\begin{abstract}
For any odd prime~$p$, the Galois group of the maximal unramified pro-$p$-extension of an imaginary quadratic field is a Schur $\sigma$-group. But Schur $\sigma$-groups can also be constructed and studied abstractly. We prove that if $p>3$, any Schur $\sigma$-group of Zassenhaus type $(3,3)$, for which every open subgroup has finite abelianization, is isomorphic to an open subgroup of a form of $\PGL_2$ over~$\BQ_p$. Combined with earlier work on an analogue of the Cohen-Lenstra heuristic for Schur $\sigma$-groups, or with the Fontaine-Mazur conjecture, this lends credence to the ``if'' part of a conjecture of McLeman.
\end{abstract}

{\renewcommand{\thefootnote}{}
\footnotetext{MSC classification: 11R11 (11R32, 20D15, 20E18, 20F05)}
%
}

\newpage
\renewcommand{\baselinestretch}{0.6}\normalsize
\tableofcontents
\renewcommand{\baselinestretch}{1.0}\normalsize

\section{Introduction}
\label{Intro}

For any odd prime $p$ and any imaginary quadratic field~$K$, the Galois group $G_K$ of the maximal unramified pro-$p$-extension of~$K$ is called the \emph{$p$-tower group} associated to~$K$. Collecting some of the known properties of $G_K$ one arrives at the notion of \emph{Schur $\sigma$-groups}, which can be constructed and studied abstractly. 

So let $G$ be a Schur $\sigma$-group with minimal number of generators~$d_G$. Following the work of Golod and Shafarevich \cite{GolodShafarevich1964} and many others, current knowledge says that $G$ is finite if $d_G\le1$, infinite if $d_G\ge3$, and in the case $d_G=2$ it is infinite unless its Zassenhaus type is (3,3) or (3,5) or (3,7) (see McLeman \cite[\S2]{McLeman2008}). Those remaining cases are not yet settled completely, but McLeman \cite[Conj.\,2.9]{McLeman2008} conjectured that any $p$-tower group of Zassenhaus type (3,3) should be finite. The present paper sheds some light on this in the case $p>3$.

\medskip
Thus we study Schur $\sigma$-groups of Zassenhaus type (3,3). To begin with, these are pro-$p$-groups with an action of a group $\{1,\sigma\}$ of order~$2$, in other words they are \emph{$\sigma$-pro-$p$-groups} according to \cite[\S4]{PinkRubio2025}. Next, any element $x$ of a $\sigma$-pro-$p$-group that satisfies ${}^\sigma x=x^{-1}$ is called \emph{odd} by \cite[\S3]{PinkRubio2025}. So let $F_2$ denote the free pro-$p$-group on two generators, turned into a $\sigma$-pro-$p$-group by making the generators odd. Then a factor group $G$ of $F_2$ described by two odd relations in the step $D_3(F_2)$ of the Zassenhaus filtration of $F_2$, whose residue classes modulo $D_4(F_2)$ are $\BF_p$-linearly independent, is called a \emph{weak Schur $\sigma$-group of Zassenhaus type (3,3).} If in addition the abelianization of $G$ is finite, then $G$ is a Schur $\sigma$-group (see \cite[Prop.\,11.2]{PinkRubio2025}).

Any $p$-tower group has the additional property that every open subgroup has finite abelianization, in other words it is a \emph{strong Schur $\sigma$-group} in the sense of \cite[Def.\,11.1]{PinkRubio2025}. We therefore concentrate on strong Schur $\sigma$-groups of Zassenhaus type $(3,3)$.

\medskip
From here on assume that $p>3$. Then our first main result states that every infinite strong Schur $\sigma$-group of Zassenhaus type $(3,3)$ is isomorphic to an open subgroup of a form of $\PGL_2$ over~$\BQ_p$ (see Proposition \ref{InfStrong33ThenA1}). 
As the Fontaine-Mazur Conjecture \cite[Conj.\,5b]{FontaineMazur1995} forbids this for $p$-tower groups, the ``if'' part of McLeman's conjecture  \cite[Conj.\,2.9]{McLeman2008} becomes a direct consequence of the Fontaine-Mazur Conjecture.

We also deduce that there exist only countably many $\sigma$-isomorphism classes of infinite strong Schur $\sigma$-group of Zassenhaus type $(3,3)$ (see Proposition \ref{InfStrSchSig33Fin}). This is significant in connection with the joint article \cite{PinkRubio2025} with Rubio, where we endow the set of $\sigma$-isomorphism classes of weak Schur $\sigma$-groups with a probability measure that serves as an analogue of the Cohen-Lenstra heuristic. In the present paper we also prove that a singleton $\{[G]\}$ in this probability space has measure $0$ if and only if $G$ is infinite (see Proposition \ref{OpenPosMeasFinite}). Using the countability result mentioned above we conclude that the set of $\sigma$-isomorphism classes of infinite weak Schur $\sigma$-groups of Zassenhaus type (3,3) is a closed subset of measure~$0$ (see Theorem \ref{I33SG4}). Assuming our heuristic to be correct, this implies that the ``if'' part of McLeman's conjecture holds on average for $p>3$.

\medskip
The paper is structured as follows. In Sections \ref{Not} and \ref{SPPG} we recall basic facts on pro-$p$-groups and $\sigma$-pro-$p$-groups, respectively. In Section \ref{SSG} we review results on weak and strong Schur $\sigma$-groups from \cite{PinkRubio2025} and prove that a singleton $\{[G]\}\subset\Sch$ is open if and only if $G$ is finite.

From Section \ref{33G} onwards we assume that $p>3$. We introduce more general pro-$p$-groups of type (3,3) by an ad hoc condition that is preserved under taking the quotient by any normal subgroup of $G$ that is contained in $D_3(G)$. In Section \ref{F33G} we study finite $p$-groups of type (3,3) and in Section \ref{P33G} we transfer these results to pro-$p$-groups of type (3,3) by taking inverse limits. We also prove that every infinite pro-$p$-group of type (3,3), for which every open subgroup has finite abelianization, is isomorphic to an open subgroup of a form of $\PGL_2$ over~$\BQ_p$. 

In Section \ref{P33SG} we deduce our main results that were described above,
and in the final Section \ref{pTG} we briefly discuss the consequences for $p$-tower groups.


%

\section{Pro-$p$-groups}
\label{Not}

Throughout this article, all homomorphisms of pro-$p$-groups are tacitly assumed to be continuous, and all subgroups are tacitly assumed to be closed. Thus by the subgroup generated by a subset of a pro-$p$-group we always mean the closure of the abstract subgroup generated by that set. In particular, when we say that a pro-$p$-group is generated by certain elements, we mean that it is topologically generated by them. 

\medskip
Consider a pro-$p$-group~$G$. For any elements $a,b\in G$ we abbreviate ${}^ab := aba^{-1}$ and $[a,b] := aba^{-1}b^{-1}$, and we write $a\sim b$ if and only if the elements are conjugate. 
For any elements $a_1,\ldots a_n\in G$ we let $\langle a_1,\ldots,a_n\rangle$ denote the subgroup generated by $a_1,\ldots a_n$.
For any subsets $A,B\subset G$ we let $[A,B]$ denote the subgroup generated by the subset $\{[a,b]\mid a\in A,\; b\in B\}.$ 
The notation $H<G$ means that $H$ is a subgroup of~$G$, but not necessarily a proper subgroup.
Likewise the notation $A\subset B$ means that $A$ is a subset of~$B$, but not necessarily a proper subset.

\medskip
We will need to use several canonically defined subgroups of~$G$. First, let $K_i(G)$ denote the descending central series defined by $K_1(G) := G$ and $K_{i+1}(G) := [G,K_i(G)]$ for all $i\ge1$ (see Berkovich \cite[page xi]{Berkovich2008}).
Second, for any $i\ge0$ let $\Agemo_i(G)$ denote the subgroup generated by the subset $\{a^{p^i}\mid a\in G\}$ (see \cite[page xiii]{Berkovich2008}).
Third, for any $i\ge1$ let $D_i(G)$ denote the $i$-th dimension subgroup of $G$ (see \cite[\S11]{DdSMS2003}). These make up the Zassenhaus filtration and are characterized by \cite[Thm.\,11.2]{DdSMS2003} as
\UseTheoremCounterForNextEquation
\begin{equation}\label{Zass}
D_i(G)\ =\ \prod_{jp^k\ge i} \Agemo_k(K_j(G)).
\end{equation}
For any letter $X\in\{K,\Agemo,D\}$ and each appropriate $i$ the group $X_i(G)$ is a characteristic subgroup of~$G$. Also the intersection $\bigcap_i X_i(G)$ is the trivial subgroup $\{1\}$. We denote the associated graded subquotients by 
$$\gr^X_i(G)\ :=\ X_i(G)/X_{i+1}(G).$$
By construction $\gr_i^K(G)$ is abelian and $\gr_i^D(G)$ is an $\BF_p$-vector space. If $G$ is finitely generated, each $D_i(G)$ is a finitely generated open subgroup and hence each $\gr_i^D(G)$ has finite dimension over~$\BF_p$. Also $D_2(G)$ is the Frattini subgroup of~$G$ and $\gr^D_1(G) = G/D_2(G)$ is the maximal quotient of $G$ that is an $\BF_p$-vector space. Moreover $d_G := \dim_{\BF_p}\gr^D_1(G)$ is the minimal number of generators of~$G$. 
\section{$\sigma$-Pro-$p$-groups}
\label{SPPG}

In this section we recall the notion and some basic properties of $\sigma$-pro-$p$ groups that were developed in \cite[\S3--4]{PinkRubio2025}, and prove a few more.

\medskip
Fix a prime $p>2$ and a finite group $\{1,\sigma\}$ of order~$2$. We call a pro-$p$-group $G$ with an action of $\{1,\sigma\}$ a \emph{$\sigma$-pro-$p$-group} and denote the action by $a\mapsto{}^\sigma a$. A $\sigma$-equivariant homomorphism between $\sigma$-pro-$p$-groups is called a \emph{$\sigma$-homomorphism}. This defines a \emph{category of $\sigma$-pro-$p$-groups} which possesses kernels and images and short exact sequences. 
The notions of \emph{$\sigma$-isomorphisms} and \emph{$\sigma$-automorphisms} are defined in the obvious way. The group of $\sigma$-automorphisms of a $\sigma$-pro-$p$-group $G$ is denoted $\Aut_\sigma(G)$.

\medskip
Consider a $\sigma$-pro-$p$-group~$G$. Then for any $\sigma$-invariant normal subgroup $N\triangleleft G$ the factor group $G/N$ is a $\sigma$-pro-$p$-group and the projection map $G\onto G/N$ is a $\sigma$-homomorphism. In particular this holds for every characteristic subgroup of $G$, for instance for the center $Z(G)$ and the commutator subgroup $[G,G]$ with the abelianization $G_\ab := G/[G,G]$. 

\medskip
To $G$ we associate the natural closed subsets
\UseTheoremCounterForNextEquation
\begin{eqnarray}\label{G+Def}
G^+\ :=\ G^{+1} &\!:=\!& \bigl\{ a\in G\bigm| {}^\sigma a = a\bigr\}, \\[3pt]
\UseTheoremCounterForNextEquation\label{G-Def}
G^-\ :=\ G^{-1} &\!:=\!& \bigl\{ a\in G\bigm| {}^\sigma a = a^{-1}\bigr\}.
\end{eqnarray}
We call the elements of $G^+$ \emph{even} and the elements of $G^-$ \emph{odd}. We call $G$ \emph{totally even} if $G=G^+$, and \emph{totally odd} if $G=G^-$. While $G^+$ is always a subgroup of~$G$, the subset $G^-$ is in general not, though it is invariant under conjugation by~$G^+$.
Both subsets are invariant under~$\sigma$, and any $\sigma$-homomorphism $G\to H$ induces a map $G^\epsilon\to H^\epsilon$ for every $\epsilon\in\{\pm1\}$. By \cite[Props.\,3.4, 4.1--2]{PinkRubio2025} 
we have:

\begin{Prop}\label{SPPGExact}
For any $\sigma$-invariant normal subgroup $N\triangleleft G$ and any $\epsilon\in\{\pm1\}$ the projection $G\onto G/N$ induces a surjective map $G^\epsilon \onto (G/N)^\epsilon$. If moreover $G$ is finite, then all fibers have the same cardinality $|N^\epsilon|$ and we have
$$|G^\epsilon|\ =\ |(G/N)^\epsilon|\cdot |N^\epsilon|.$$
\end{Prop}

\begin{Prop}\label{SPPGGrad}
The product map induces homeomorphisms
$$\begin{array}{c}
G^+\times G^- \longto G,\ (a,b)\mapsto ab. \\[3pt]
G^-\times G^+ \longto G,\ (a,b)\mapsto ab.
\end{array}$$
\end{Prop}

\medskip
We will also need the following facts:

\begin{Prop}\label{TotOddAbel}
Any totally odd $\sigma$-pro-$p$-group $G$ is abelian.
\end{Prop}

\begin{Proof}
By assumption the map $G\to G$, $a\mapsto a^{-1}$ is equal to $\sigma$ and hence a homomorphism. It is well-known that this implies that $G$ is abelian.
\end{Proof}

%


\begin{Prop}\label{StrongAutFinThenFin}
Let $G$ be a $\sigma$-pro-$p$-group such that $\Aut_\sigma(G)$ is finite and for every open subgroup $H<G$ the abelianization $H_\ab$ is finite. Then $G$ is finite.
\end{Prop}

\begin{Proof}
For any $a\in G^+$ the associated inner automorphism of $G$ commutes with~$\sigma$, yielding a natural homomorphism $G^+\to\Aut_\sigma(G)$. By definition the kernel of this homomorphism is the even part $Z(G)^+$ of the center $Z(G)$ of~$G$. As $\Aut_\sigma(G)$ is finite, that has finite index in~$G^+$ and is therefore open in~$G^+$. Since $G^+$ carries the subspace topology from~$G$, there thus exists an open normal subgroup $N\triangleleft G$ that satisfies $N\cap G^+ \subset Z(G)^+$. After replacing $N$ by $N\cap{}^\sigma N$ we may in addition assume that $N$ is $\sigma$-invariant.

Setting $M:=N\cap Z(G)$ this means that $N^+ = M^+$. By Proposition \ref{SPPGExact} the factor group $N/M$ thus satisfies $(N/M)^+=1$. By Proposition \ref{SPPGGrad} it is therefore totally odd, and so it is abelian by Proposition \ref{TotOddAbel}. Thus $N/M$ is a factor group of~$N_\ab$. Since $N_\ab$ is finite by the second assumption on~$G$, it follows that $N/M$ is finite. Therefore $M$ is open in $N$ and hence also in~$G$. As $M$ is abelian, using the second assumption on $G$ again we deduce that $M$ is finite. Thus $G$ is finite, as desired.
\end{Proof}

\medskip
Finally, for every integer $n\ge0$ we let $F_n$ denote the free pro-$p$-group on $n$ generators, turned into a $\sigma$-pro-$p$-group by making the generators odd.

\section{Schur $\sigma$-groups}
\label{SSG}

In this section we recall some material on weak and strong Schur $\sigma$-groups from \cite[\S7--11]{PinkRubio2025} and prove a further result concerning measures.

\medskip
First, a $\sigma$-pro-$p$-group $G$ with the properties
\UseTheoremCounterForNextEquation
\begin{equation}\label{SchConsDef}
\scriptstyle\biggl\{\displaystyle
\begin{array}{c}
\hbox{$\dim_{\BF_p}\!H^2(G,\BF_p) \le \dim_{\BF_p}\!H^1(G,\BF_p) < \infty$} \\[3pt]
\hbox{and $\sigma$ acts by $-1$ on both spaces}
\end{array}\scriptstyle\biggr\}
\end{equation}
is called a \emph{weak Schur $\sigma$-group}. We construct a certain set $\{ F_n/N_\rtuple \mid n\ge0,\ \rtuple\in(F_n^-)^n \}$ of Schur $\sigma$-groups such that each Schur $\sigma$-group is $\sigma$-isomorphic to one in this set. The set $\Sch$ (pronounce \emph{Schur}) of $\sigma$-isomorphism classes in this set can thus be viewed as the set of all weak Schur $\sigma$-groups up to $\sigma$-isomorphism, avoiding set-theoretic complications. See \cite[Def.\,7.8]{PinkRubio2025}.

\medskip
Next, for any $i\ge2$ and any weak Schur $\sigma$-group $G$ consider the subset
\UseTheoremCounterForNextEquation
\begin{equation}\label{UnGDef}
U_{i,G}\ :=\ \bigl\{ [H]\in\Sch \bigm| H/D_i(H) \cong G/D_i(G) \bigr\}.
\end{equation}
We show that these sets form a basis for a unique topology on~$\Sch$, which is Hausdorff and totally disconnected, and for which each set $U_{i,G}$ is open and closed. See \cite[Prop.\,8.2]{PinkRubio2025}. 
In \cite[Prop.\,8.7 (a)]{PinkRubio2025} 
we prove:

\begin{Prop}\label{FinSchurOpenClosed}
If $G$ is finite, the subset $\{[G]\} \subset \Sch$ is open and closed.
\end{Prop}

\medskip
Next, in \cite[\S9]{PinkRubio2025} we construct a certain probability measure $\mu_\infty$ on the Borel $\sigma$-algebra of $\Sch$. Also, for any integer $k\ge0$ as well as for $k=\infty$ consider the positive real number 
\UseTheoremCounterForNextEquation
\begin{equation}\label{CkDef}
C_k\ :=\ \prod_{i=1}^k (1-p^{-i})\ >\ 0.
\end{equation}
In \cite[Prop.\,10.8]{PinkRubio2025} 
we prove:

\begin{Prop}\label{UiGMeasure}
Let $n$ be the minimal number of generators of $G$ and $m$ the minimal number of relations of $G/D_i(G)$ as a factor group of $F_n/D_i(F_n)$. Then
$$\mu_\infty(U_{i,G})\ =\ \frac{C_\infty}{C_{n-m}}\cdot \frac{1}{\bigl|\Aut_\sigma(G/D_i(G))\bigr|}.$$
\end{Prop}

Next, a weak Schur $\sigma$-group $G$ such that for every open subgroup $H<G$ the abelianization $H_\ab$ is finite is called a \emph{strong Schur $\sigma$-group}. 
(If one only requires that $G_\ab$ is finite one gets a \emph{Schur $\sigma$-group} by 
\cite[Prop.\,11.2]{PinkRubio2025}.)
Clearly every finite weak Schur $\sigma$-group is a strong Schur $\sigma$-group. See \cite[Def.\,11.1]{PinkRubio2025}. 
In \cite[Thm.\,11.4]{PinkRubio2025} 
we prove:

\begin{Thm}\label{StrongThm}
The subset $\Sch^\strong$ of $\sigma$-isomorphisms classes of strong Schur $\sigma$-groups has measure~$1$.
\end{Thm}

Now we can add to this one further result:

\begin{Prop}\label{OpenPosMeasFinite}
For any weak Schur $\sigma$-group~$G$ the following are equivalent:
\begin{enumerate}
\item[(a)] The set $\{[G]\}$ is open in $\Sch$.
\item[(b)] The set $\{[G]\}$ has measure $>0$.
\item[(c)] The group $G$ is finite.
\end{enumerate}
\end{Prop}

\begin{Proof}
The implication (c)$\Rightarrow$(a) is contained in Proposition \ref{FinSchurOpenClosed}. 

Next, since $\Sch$ is Hausdorff, the subset $\{[G]\}$ is always closed and hence measurable. If it is open, by the definition of the topology on $\Sch$ it is equal to $U_{i,G}$ for some $i>0$.
Proposition \ref{UiGMeasure} thus implies that its measure is~$>0$, proving the implication (a)$\Rightarrow$(b).

Now assume that $\mu_\infty(\{[G]\})>0$. Since the subset of $\sigma$-isomorphisms classes of weak Schur $\sigma$-groups that are not strong has measure~$0$ by Theorem \ref{StrongThm}, the assumption implies that $G$ is a strong Schur $\sigma$-group.
On the other hand, for every $i\ge2$ Proposition \ref{UiGMeasure} implies that 
$$\bigl|\Aut_\sigma(G/D_i(G))\bigr|
\ =\ \frac{C_\infty}{C_{n-m}}\cdot\frac{1}{\mu_\infty(\{[G]\})}$$
for some $0\le m\le n$. As $n=\dim_{\BF_p}G/D_2(G)$ is independent of $i\ge2$, there are only finitely many possible values for $C_{n-m}$ and they are all $>0$. It follows that $|\Aut_\sigma(G/D_i(G))|$ is bounded independently of~$i$, say by a number~$b$. 
Since $\Aut_{\sigma}(G)$ is the filtered inverse limit of the groups $\Aut_{\sigma}(G/D_i(G))$, it follows that $\Aut_\sigma(G)$ is finite of cardinality $\le b$.
By Proposition \ref{StrongAutFinThenFin} the group $G$ is therefore finite, proving the implication (b)$\Rightarrow$(c).
\end{Proof}

\begin{Cor}\label{CountInfNull}
Any at most countable subset of $\{[G]\in\Sch\mid G\ \hbox{is infinite}\}$ has measure~$0$.
\end{Cor}

\begin{Proof}
By Proposition \ref{OpenPosMeasFinite} the set $\{[G]\}$ has measure $0$ whenever $G$ is infinite, and any countable union of sets of measure $0$ has measure~$0$.
\end{Proof}

\begin{Rem}\label{CardRem1}\rm
For any $n\ge1$ the set $F_n^-$ has the cardinality of the continuum. Indeed, being the inverse limit of the finite sets $(F_n/D_i(F_n))^-$ over all $i\ge1$, it has at most that cardinality, and conversely, since it surjects to $F_{n,\ab} \cong \BZ_p^n$, it has at least that cardinality. 
As by construction we have a surjective map $\bigsqcup_{n\ge0}(F_n^-)^n \onto \Sch$, it follows that $\Sch$ has at most the cardinality of the continuum.
\end{Rem}

\begin{Rem}\label{CardRem2}\rm
As there are only finitely many isomorphism classes of finite $\sigma$-$p$-groups of any given order, the subset $\{[G]\in\Sch\mid G\ \hbox{is finite}\}$ is countable. 
On the other hand, by Proposition \ref{CountInfNull} any subset of measure $>0$ of $\{[G]\in\Sch\mid G\ \hbox{is infinite}\}$ is uncountable. Since any weak Schur $\sigma$-group $G$ with minimal number of generators $\ge3$ is infinite, it follows that $\Sch^\strong$ is uncountable. 
\end{Rem}

\begin{Rem}\label{CardRem3}\rm
One can prove that each of $\Sch^\strong$ and $\Sch\setminus\Sch^\strong$ has the cardinality of the continuum. As this seems only mildly interesting, we leave the details as an exercise.
\end{Rem}

\section{Groups of type (3,3)}
\label{33G}

First consider a weak Schur $\sigma$-group~$G$. Recall (for instance from Koch-Venkov \cite{KochVenkov1974}) that $G$ has Zassenhaus type (3,3) if and only if it can be described by two generators modulo two relations, which lie in the filtration step $D_3(F_2)$ of the free pro-$p$-group~$F_2$ and whose residue classes in the subquotient $\gr^D_3(F_2)$ are $\BF_p$-linearly independent.
To describe the precise impact of this condition, recall that, since $p>2$, by \cite[Thm.\,5]{McLeman2009} we have 
\UseTheoremCounterForNextEquation
\begin{equation}\label{F2Grad}
\dim_{\BF_p}\gr^D_i(F_2)\ =\ 
\scriptstyle\left\{\kern-2pt\textstyle 
\begin{array}{ll}
2 &\hbox{if $i=1$,}\\[3pt]
1 &\hbox{if $i=2$,}\\[3pt]
2 &\hbox{if $i=3$ and $p>3$,}\\[3pt]
4 &\hbox{if $i=3$ and $p=3$.}
\end{array}\right.
\end{equation}
For $p=3$ there are therefore many more subtypes of the Zassenhaus type (3,3) than for $p>3$. As several other key arguments in Section \ref{F33G} also fail for $p=3$, throughout the rest of this paper we impose the assumption:

\begin{Ass}\label{PGreaterThanThree}
We have $p>3$.
\end{Ass}

{}From \eqref{F2Grad} we deduce that any weak Schur $\sigma$-group $G$ of Zassenhaus type (3,3) satisfies
\UseTheoremCounterForNextEquation
\begin{equation}\label{33Def}
\dim_{\BF_p}\gr^D_i(G)\ =\ 
\scriptstyle\left\{\kern-2pt\textstyle 
\begin{array}{ll}
2 &\hbox{if $i=1$,}\\[3pt]
1 &\hbox{if $i=2$,}\\[3pt]
0 &\hbox{if $i=3$.}
\end{array}\right.
\end{equation}
We will analyze such a group by first studying all its finite quotients. Note that the condition \eqref{33Def} is preserved under taking the quotient by any normal subgroup of $G$ that is contained in $D_3(G)$. For convenience we therefore give this condition a name, and call any pro-$p$-group satisfying \eqref{33Def} a \emph{pro-$p$-group of type (3,3)}.

\medskip
So consider an arbitrary pro-$p$-group $G$ of type (3,3). The condition \eqref{33Def} for $i=1$ means that the minimal number of generators of $G$ is~$2$. Writing $G$ as a quotient of the free pro-$p$-group~$F_2$, the other conditions mean that all relations lie in $D_3(F_2)$ and their residue classes generate $\gr^D_3(F_n)$. At this point we do not care about whether there are further relations in lower steps of the Zassenhaus filtration or not.

The conditions imply that $G/D_3(G)$ is a non-abelian group of order~$p^3$ and exponent~$p$. As there is a unique isomorphism class of such groups, we deduce that
\UseTheoremCounterForNextEquation
\begin{equation}\label{33GIsom}
G/D_3(G)\ \cong\ \left(\begin{smallmatrix}
1 & * & * \\[2pt] 0 & 1 & * \\[2pt] 0 & 0 & 1
\end{smallmatrix}\right)
\ <\ \GL_3(\BF_p).
\end{equation}
For the following we pick generators $x$ and $y$ of~$G$, whose residue classes therefore form a basis of $\gr^D_1(G)$ over~$\BF_p$. We also fix an element $z\in D_2(G)\setminus D_3(G)$, whose residue class therefore forms a basis of $\gr^D_2(G)$ over~$\BF_p$. For example we might take $z=[x,y]$, but we prefer to leave ourselves more freedom by allowing an arbitrary~$z$.

\section{Finite $p$-groups of type (3,3)}
\label{F33G}

In this section we fix a finite $p$-group $G$ of type (3,3).

\begin{Prop}\label{F33G1}
We have 
$$D_3(G) = \Agemo_1(G).$$
\end{Prop}

\begin{Proof}
Since $p>2$, the formula \eqref{Zass} implies that $D_3(G)=K_3(G)\cdot \Agemo_1(G)$. The desired equation is therefore equivalent to $K_3(G)\subset \Agemo_1(G)$.
To prove this we may replace $G$ by $G/\Agemo_1(G)$, which is again a finite $p$-group of type (3,3). Thus we may assume that $\Agemo_1(G)=1$ and need to prove that $K_3(G)=1$.

With \eqref{Zass} the assumption now implies that $D_i(G)=K_i(G)$ for all $i\ge1$. The condition for $i=3$ in \eqref{33Def} thus means that $K_3(G)=K_4(G)$. Next, if $K_i(G)=K_{i+1}(G)$ for some integer $i\ge3$, the computation
$$K_{i+1}(G)\ =\ [G,K_i(G)]\ =\ [G,K_{i+1}(G)]\ =\ K_{i+2}(G)$$
implies that $K_{i+1}(G)=K_{i+2}(G)$ as well. By induction we therefore deduce that $K_3(G)=K_i(G)$ for all $i\ge3$. But since $G$ is a finite $p$-group, we have $K_i(G)=1$ for some $i\gg0$. Thus $K_3(G)=1$, and we are done.
\end{Proof}

\begin{Prop}\label{F33G2}
The group $G$ is a \emph{regular $p$-group} in the sense of Berkovich \cite[\S7 Def.\,1]{Berkovich2008}, that is, for any $a,b\in G$ we have
$$a^pb^p \equiv (ab)^p \bmod \Agemo_1(K_2(\langle a,b\rangle)).$$
\end{Prop}

\begin{Proof}
Since $p>3$, the isomorphism \eqref{33GIsom} and Proposition \ref{F33G1} imply that $[G:\Agemo_1(G)] = {p^3<p^p}$. Thus $G$ is absolutely regular in the sense of \cite[\S9 Def.\,1]{Berkovich2008} and hence regular by Hall's regularity criterion  \cite[Thm.\,9.8 (a)]{Berkovich2008}.
\end{Proof}

\begin{Prop}\label{F33G3}
For any subgroup $H<G$ and any $i,j\ge0$ we have 
$$\Agemo_i(\Agemo_j(H)) = \Agemo_{i+j}(H).$$
\end{Prop}

\begin{Proof}
By the definition of regularity any subgroup $H<G$ is itself regular (see also \cite[Thm.\,7.1 (a)]{Berkovich2008}). For any $i\ge0$ we therefore have $\Agemo_i(H) = \{h^{p^i}\mid h\in H\}$ by \cite[Thm.\,7.2~(c)]{Berkovich2008}). This directly implies the claim.
\end{Proof}


\begin{Prop}\label{F33G5}
For any $a,b\in G$ we have
\begin{enumerate}
\item[(a)] $a^pb^p \equiv (ab)^p \bmod \Agemo_1(D_2(G))$, and
\item[(b)] $a^pb^p \equiv (ab)^p \bmod \Agemo_2(G)$ if $a\in D_2(G)$ or $b\in D_2(G)$.
\end{enumerate}
\end{Prop}

\begin{Proof}
Since $K_2(\langle a,b\rangle) \subset K_2(G) \subset D_2(G)$, the first congruence follows directly from Proposition \ref{F33G2}. If in addition $a\in D_2(G)$ or $b\in D_2(G)$, the fact that $D_2(G)/D_3(G)$ is the center of $G/D_3(G)$ implies that $K_2(\langle a,b\rangle)\subset D_3(G)$. Using Propositions \ref{F33G1} and \ref{F33G3} it follows that
$$\Agemo_1(K_2(\langle a,b\rangle))\ \subset\ \Agemo_1(D_3(G))\ =\ \Agemo_1(\Agemo_1(G))\ =\ \Agemo_2(G).$$
Thus the second congruence again follows directly from Proposition \ref{F33G2}. 
\end{Proof}

\begin{Lem}\label{F33G6}
We have 
$$\Agemo_1(D_2(G))\ \subset\ \langle z^p\rangle\cdot \Agemo_2(G).$$
\end{Lem}

\begin{Proof}
Since $D_2(G)/D_3(G)$ is generated by the residue class of~$z$, any element of $D_2(G)$ can be written in the form $a=z^ib$ for some $i\in\BZ$ and some $b\in D_3(G)$. Proposition \ref{F33G5} (b) then implies that $a^p \equiv z^{ip}b^p \bmod \Agemo_2(G)$. Here $b^p$ lies in $\Agemo_1(D_3(G)) = \Agemo_1(\Agemo_1(G)) = \Agemo_2(G)$ by Propositions \ref{F33G1} and \ref{F33G3}; hence $a^p$ lies in the subgroup $\langle z^p\rangle\cdot \Agemo_2(G)$. Varying $a$ it follows that the subgroup $\Agemo_1(D_2(G))$ is contained in $\langle z^p\rangle\cdot \Agemo_2(G)$, as desired.
\end{Proof}

\begin{Prop}\label{F33G7}
The subgroup $\Agemo_1(G)$ is generated by $x^p,y^p,z^p$.
\end{Prop}

\begin{Proof}
Since $G/D_2(G)$ is abelian and generated by the residue classes of~$x$ and $y$, any element of $G$ can be written in the form $a=x^iy^jb$ for some $i,j\in\BZ$ and some $b\in D_2(G)$. Proposition \ref{F33G5} (a) then implies that $a^p \equiv x^{ip}y^{jp}b^p \bmod \Agemo_1(D_2(G))$. Since $b^p\in \Agemo_1(D_2(G))$, this means that $a^p \equiv x^{ip}y^{jp} \bmod \Agemo_1(D_2(G))$. By Lemma \ref{F33G6} this shows that $a^p$ lies in the subgroup $\langle x^p,y^p,z^p\rangle\cdot \Agemo_2(G)$. Varying $a$ it follows that the subgroup $\Agemo_1(G)$ is contained in $\langle x^p,y^p,z^p\rangle\cdot \Agemo_2(G)$. Conversely the latter group is contained in $\Agemo_1(G)$ by the definition of the subgroups $\Agemo_i(G)$. Combining this with Proposition~\ref{F33G3} we find that
$$\Agemo_1(G)\ =\ \langle x^p,y^p,z^p\rangle\cdot \Agemo_2(G)
\ =\ \langle x^p,y^p,z^p\rangle\cdot \Agemo_1(\Agemo_1(G)).$$
Since $\Agemo_1(\Agemo_1(G))$ is contained in the Frattini subgroup of $\Agemo_1(G)$, it follows that $\Agemo_1(G) = \langle x^p,y^p,z^p\rangle$, as desired.
\end{Proof}


\begin{Prop}\label{F33G8}
The group $\Agemo_1(G)$ is \emph{powerful} in the sense of \cite[\S26]{Berkovich2008}, that is, we have
$$K_2(\Agemo_1(G)) \subset \Agemo_1(\Agemo_1(G)).$$
\end{Prop}

\begin{Proof}
This follows from regularity by \cite[\S26 Exercise 2]{Berkovich2008} or \cite[page 485 (2)]{LubotzkyMann1987a}.
\end{Proof}

\begin{Prop}\label{F33G9a}
Any element of $\Agemo_1(G)$ can be written in the form $x^{pi}y^{pj}z^{pk}$ with ${i,j,k\in\BZ}$.
\end{Prop}

\begin{Proof}
Direct combination of Propositions \ref{F33G7} and \ref{F33G8} and \cite[Cor.\,26.12]{Berkovich2008}.
\end{Proof}

\begin{Prop}\label{F33G9}
Any element of $G$ can be written in the form $x^iy^jz^k$ with ${i,j,k\in\BZ}$.
\end{Prop}

\begin{Proof}
Take any $a\in G$. By the description of $G/D_3(G)$ from \eqref{33GIsom}, we can choose integers $i_0,j_0,k_0$ such that $a$ is congruent to $b:= x^{i_0}y^{j_0}z^{k_0}$ modulo $D_3(G)$. It then suffices to find $i,j,k\in\BZ$ such that 
$$a\ =\ x^{i_0+pi}y^{j_0+pj}z^{k_0+pk}.$$
This equation is equivalent to
\begin{eqnarray*}
a'\ :=\ b^{-1}a
&=& (z^{-k_0}y^{-j_0}x^{-i_0}) (x^{i_0+pi}y^{j_0+pj}z^{k_0+pk})\\
&=& z^{-k_0}y^{-j_0}x^{pi}y^{j_0}y^{pj}z^{k_0}z^{pk}\\
&=& z^{-k_0}y^{-j_0}x^{pi}y^{j_0}z^{k_0}\cdot z^{-k_0}y^{pj}z^{k_0}\cdot z^{pk}\\
&=& (x')^{pi}(y')^{pj}z^{pk}
\end{eqnarray*}
with $x' := z^{-k_0}y^{-j_0}xy^{j_0}z^{k_0}$ and $y' := z^{-k_0}yz^{k_0}$. Here the residue classes of $x'$ and $y'$ again generate the abelian group $G/D_2(G)$; hence all the above results hold equally for $(x',y',z)$ in place of $(x,y,z)$. On the other hand we have $a'=b^{-1}a\in D_3(G)=\Agemo_1(G)$ by Proposition \ref{F33G1}. By Proposition \ref{F33G9a} we can therefore write $a'=(x')^{pi}(y')^{pj}z^{pk}$ for some $i,j,k\in\BZ$, as desired.
\end{Proof}

\begin{Prop}\label{F33G10}
For every $k\ge1$ the group $\gr^\Agemo_k(G)$ is an $\BF_p$-vector space of dimension $\le3$ generated by the residue classes of $x^{p^k},y^{p^k},z^{p^k}$ and the map $a\mapsto a^p$ induces a surjective homomorphism $\gr_k^\Agemo(G)\onto\gr_{k+1}^\Agemo(G)$.
\end{Prop}

\begin{Proof}
Since $H := \Agemo_1(G)$ is powerful, \cite[Thm.\,26.12 (b)]{Berkovich2008} implies that the subgroups $P_k(H)$ defined just before \cite[Lem.\,26.8]{Berkovich2008} satisfy $P_k(H) = \Agemo_{k-1}(P_1(H)) = \Agemo_{k-1}(H)$. By Proposition \ref{F33G3} this means that $P_k(H) = \Agemo_k(G)$. By \cite[Lem.\,26.8 (a)]{Berkovich2008} it follows that $\gr^\Agemo_k(G)$ is an $\BF_p$-vector space and the map $a\mapsto a^p$ induces a surjective homomorphism $\gr_k^\Agemo(G)\onto\gr_{k+1}^\Agemo(G)$. By Proposition \ref{F33G7} the group $\gr_1^\Agemo(G)$ is generated by the residue classes of $x^p,y^p,z^p$ and therefore of dimension $\le3$. The corresponding statements then follow for all~$k$.
\end{Proof}

\section{Pro-$p$-groups of type (3,3)}
\label{P33G}

Now consider a pro-$p$-group $G$ of type (3,3) with elements $x,y,z$ as in Section~\ref{33G}. For any open normal subgroup $H\triangleleft G$ that is contained in $D_3(G)$, the factor group $G/H$ is a finite $p$-group of type (3,3) and the images of $x,y,z$ in $G/H$ satisfy the corresponding properties. By taking inverse limits Propositions \ref{F33G1} and \ref{F33G9} and \ref{F33G10} imply:

\begin{Prop}\label{P33G1}
We have $D_3(G) = \Agemo_1(G)$.
\end{Prop}

\begin{Prop}\label{P33G9}
Any element of $G$ can be written in the form $x^iy^jz^k$ with ${i,j,k\in\BZ_p}$.
\end{Prop}

\begin{Prop}\label{P33G10}
For every $k\ge1$ the group $\gr^\Agemo_k(G)$ is an $\BF_p$-vector space of dimension $\le3$ generated by the residue classes of $x^{p^k},y^{p^k},z^{p^k}$ and the map $a\mapsto a^p$ induces a surjective homomorphism $\gr_k^\Agemo(G)\onto\gr_{k+1}^\Agemo(G)$.
\end{Prop}



{}From Proposition \ref{P33G10} we deduce that there exists a unique integer $d\le3$ and an integer $k_0\ge1$ such that for every $k\ge k_0$ we have $\dim_{\BF_p}\gr^\Agemo_k(G)=d$ and $\gr_k^\Agemo(G)\isoto\gr_{k+1}^\Agemo(G)$. Thus $\Agemo_{k_0}(G)$ is a group of rank~$d$ in the sense of Lazard \cite[Ch.\,III Def.\,2.1.3]{Lazard1965} and is therefore a $p$-adic analytic group of dimension $d$ by \cite[Ch.\,III Lemme\,3.1.5, Th.\,3.1.7.1]{Lazard1965}. The same conclusion then also holds for~$G$. This means that $G$ has a unique structure of a $p$-adic analytic variety (in the sense of \cite[Ch.\,III 1.3.3]{Lazard1965}) such that the group operation is given by analytic functions.
In particular $G$ possesses an associated Lie algebra $\Lie G$, which is a $\BQ_p$-vector space of dimension~$d$ and comes with a natural adjoint representation
\UseTheoremCounterForNextEquation
\begin{equation}\label{AdRep}
\Ad\colon G\longto \Aut_{\BQ_p}(\Lie G) \cong \GL_d(\BQ_p).
\end{equation}
Let $\CG\subset\GL_{d,\BQ_p}$ denote the Zariski closure of the image $\Ad(G)$.

\bigskip
For the rest of this section we impose the following additional assumption:

\begin{Ass}\label{InfStrongAss}
The group $G$ is infinite, and every open subgroup $H$ has finite abelianization~$H_\ab$.
\end{Ass}

\begin{Prop}\label{P33G1a}
Then $\Ad(G)$ is open in $\CG(\BQ_p)$ and $\CG$ is a form of $\PGL_2$ over~$\BQ_p$.
\end{Prop}

\begin{Proof}
By \cite[Cor.\,4.19]{DdSMS2003} the kernel $\ker(\Ad)$ contains an abelian subgroup of finite index in $\ker(\Ad)$. Thus if $\CG$ is finite, that is an abelian subgroup of finite index in~$G$. This subgroup is then open and therefore finite by the second part of Assumption \ref{InfStrongAss}. But this contradicts the assumption that $G$ is infinite. Therefore $\CG$ is not finite.

Its identity component $\CG^\circ$ is thus a non-trivial connected linear algebraic group over~$\BQ_p$. By construction $G^\circ := \Ad^{-1}(\CG^\circ(\BQ_p))$ is an open normal subgroup of $G$ whose image $\Ad(G^\circ)$ is Zariski dense in~$\CG^\circ$.

If $\CG^\circ$ has an abelian quotient $\CH$ of dimension $>0$, the image of $G^\circ$ in $\CH(\BQ_p)$ is Zariski dense and therefore infinite. Being also abelian, it follows that $(G^\circ)_\ab$ is infinite, again contradicting Assumption \ref{InfStrongAss}. Thus $\CG^\circ$ possesses no abelian quotient of dimension $>0$ and is therefore equal to its commutator subgroup. In particular $\CG^\circ$ is not solvable.

Being a linear algebraic group over a field of characteristic zero, its Lie algebra then satisfies $[\Lie\CG^\circ,\Lie\CG^\circ]=\Lie\CG^\circ$. On the other hand $\Ad(G^\circ)$ is a closed subgroup of $\CG^\circ(\BQ_p)$ and therefore a $p$-adic analytic subgroup with a Lie algebra $\Lie\Ad(G^\circ) \subset \Lie\CG^\circ$. Since $\Ad(G^\circ)$ is Zariski dense in $\CG^\circ(\BQ_p)$, by \cite[Ch.\,II Cor.\,7.9]{BorelLAG} it follows that $\Lie\Ad(G^\circ) = \Lie\CG^\circ$. This means that $\Ad(G^\circ)$ is open in $\CG^\circ(\BQ_p)$.

Since $\Ad$ is a $p$-adic analytic homomorphism, it follows that $\dim(\CG^\circ) = \dim(\Ad(G^\circ)) \le \dim(G^\circ) = \dim(G)\le d\le 3$. As $\CG^\circ$ is a non-solvable connected linear algebraic group, for dimension reasons it must therefore be semisimple with a root system of type~$A_1$. In particular we deduce that $\dim(\CG^\circ)=d=3$ and that $\Ad$ induces an isomorphism $\Lie(G^\circ)\isoto\Lie(\CG^\circ)$. Since $G$ acts on its Lie algebra through Lie algebra automorphisms, and the automorphism group of any Lie algebra of type $A_1$ is a form of $\PGL_2$, the group $\CG$ must be equal to that form of $\PGL_2$. In particular it follows that $\CG=\CG^\circ$ and hence $G=G^\circ$, and we are done.
\end{Proof}

\begin{Prop}\label{P33G1b}
The map $\BZ_p^3\longto G$, $(i,j,k)\mapsto x^iy^jz^k$ is bijective.
\end{Prop}

\begin{Proof}
By the preceding proof we already know that $d=3$. Thus in Proposition \ref{P33G10} the residue classes of $x^{p^k},y^{p^k},z^{p^k}$ must form a basis of $\gr^\Agemo_k(G)$ for every $k\ge1$. By induction it follows that the map
$$(\BZ/p^n\BZ)^3 \longto G/\Agemo_n(G),\ 
([i],[j],[k]) \mapsto [x^iy^jz^k]$$
is bijective for every $n\ge1$. The proposition follows by taking the inverse limit over~$n$.
\end{Proof}

\begin{Prop}\label{P33G1c}
The homomorphism $\Ad\colon G\to\CG(\BQ_p)$ is an open embedding.
\end{Prop}

\begin{Proof}
Since $\Ad\colon G\onto\Ad(G)$ is a surjective homomorphism of $p$-adic analytic groups of dimension $d=3$, its kernel must be finite. 

Suppose that this kernel is non-trivial and choose an element $g$ of order~$p$. If this element lies in $G\setminus D_2(G)$, then in all the above we can replace $x$ or $y$ by~$g$. But then Proposition \ref{P33G1b} becomes false, so this case cannot occur. Similarly, if this element lies in $D_2(G)\setminus D_3(G)$, then all the above holds with $z=g$, again contradicting Proposition \ref{P33G1b}. Thus $g$ lies in $D_3(G)$ and hence in $\Agemo_1(G)$ by Proposition~\ref{P33G1}. Since $\bigcap_k\Agemo_k(G)=\{1\}$, there then exists some $k\ge1$ with $g\in\Agemo_k(G)\setminus\Agemo_{k+1}(G)$. Thus $g$ represents a non-zero element of $\gr^\Agemo_k(G)$. But we have already seen that $\gr^\Agemo_k(G)$ and $\gr^\Agemo_{k+1}(G)$ both have dimension~$3$. The surjective homomorphism $\gr_k^\Agemo(G)\onto\gr_{k+1}^\Agemo(G)$, $[a]\mapsto[a^p]$ from Proposition \ref{P33G10} is therefore an isomorphism. Thus $g^p$ must represent a non-zero element of $\gr_{k+1}^\Agemo(G)$, contradicting the fact that $g^p=1$. 

Thus the kernel is trivial, and so $\Ad\colon G\to\CG(\BQ_p)$ is injective. As we already know that $\Ad(G)$ is open in $\CG(\BQ_p)$, we are done.
\end{Proof}

\section{Infinite Schur $\sigma$-groups of Zassenhaus type (3,3)}
\label{P33SG}

Keeping the assumption $p>3$, we now return to weak Schur $\sigma$-groups.

\begin{Prop}\label{InfStrong33ThenA1}
Every infinite strong Schur $\sigma$-group $G$ of Zassenhaus type $(3,3)$ is isomorphic to an open subgroup of a form $\CG$ of $\PGL_2$ over~$\BQ_p$, such that the action of $\sigma$ on $G$ is induced by an automorphism of order $2$ of~$\CG$.
\end{Prop}

\begin{Proof}
The conditions imply \eqref{33Def} and Assumption \ref{InfStrongAss}, so the first statement follows from Propositions \ref{P33G1a} and~\ref{P33G1c}. 
For the second observe that the action of $\sigma$ on $G$ induces an action on its Lie algebra, under which the embedding
\eqref{AdRep} is equivariant. By the construction of $\CG$ as the Zariski closure of the image it therefore also inherits an action of~$\sigma$. Since $\sigma$ acts non-trivially on the generators of~$G$, it must also act non-trivially on~$\CG$.
\end{Proof}


\begin{Prop}\label{PGL2FormSigmaFin}
Up to isomorphism there exist only finitely many forms of $\PGL_2$ over~$\BQ_p$ together with an automorphism of order~$2$.
\end{Prop}

\begin{Proof}
Every form $\CG$ of $\PGL_2$ over any field $K$ arises from a quaternion algebra $Q$ over~$K$ such that $\CG(K) \cong Q^\times/K^\times$. By the Skolem-Noether theorem, conjugation by $Q^\times$ then induces an isomorphism $Q^\times/K^\times \isoto \Aut_K(Q)$. Moreover, the fact that every automorphism of any form of $\PGL_2$ is inner yields a natural isomorphism $\CG(K) \isoto \Aut_K(\CG)$. Thus all automorphisms of order $2$ of $\CG$ comes from automorphisms of order $2$ of $Q$ over~$K$. 

To describe these assume now that $K$ has characteristic different from~$2$. Then over an algebraic closure $\bar K$ of~$K$, every element of order~$2$ of $\PGL_2(\bar K)$ is conjugate to the matrix
$\binom{\kern2pt 1\kern9pt 0\kern2pt}{0\ -1}$. Its centralizer in the matrix ring $\Mat_{2\times2}(\bar K)$ is thus a semisimple commutative $\bar K$-algebra isomorphic to $\bar K\times \bar K$. Since $Q$ is a form of that matrix ring over~$K$, it follows that the centralizer of any automorphism of order $2$ is a semisimple commutative $K$-subalgebra $L$ of $Q$ of dimension $2$ over~$K$. Thus $L$ is either isomorphic to $K\times K$ or a field extension of degree $2$ of~$K$. Moreover, the automorphism of order $2$ must then be the identity on $L$ and multiplication by $-1$ on its orthogonal complement for the reduced trace form on~$Q$. Thus the automorphism is completely determined by the subalgebra~$L$.

For the local field $K=\BQ_p$ there are precisely two isomorphism classes for the quaternion algebra~$Q$, for instance by \cite[Thm.\,13.1.6]{Voight2021}. Next $\BQ_p$ possesses only finitely many field extensions of degree $2$ up to isomorphism; 
hence there are only finitely many isomorphism classes for~$L$. 
Moreover, by the Skolem-Noether theorem in the form of \cite[Thm.\,7.7.1]{Voight2021}, any two embeddings $L\into Q$ over~$K$ are conjugate under~$Q^\times$. Thus there are only finitely many pairs $(Q,L)$ up to isomorphism, and the proposition follows.
\end{Proof}

\begin{Prop}\label{InfStrSchSig33Fin}
There exist at most countably many $\sigma$-isomorphism classes of infinite strong Schur $\sigma$-group of Zassenhaus type $(3,3)$.
\end{Prop}

\begin{Proof}
By Proposition \ref{InfStrong33ThenA1} each such group is $\sigma$-isomorphic to a $\sigma$-invariant open subgroup of a form $\CG$ of $\PGL_2$ over~$\BQ_p$ together with an automorphism of order~$2$. By Proposition \ref{PGL2FormSigmaFin} there are only finitely many possibilities for the latter up to isomorphism. Finally, for each such $\CG$ the topology of $\CG(\BQ_p)$ is second countable; hence there exist only countably many open subgroups.
\end{Proof}

\begin{Rem}\label{I33SG2}\rm
Though immaterial for our purposes, there are indeed countably many $\sigma$-isomorphism classes. More precisely, for any form $\CG$ of $\PGL_2$ over~$\BQ_p$ with an automorphism $\sigma$ of order~$2$ and any $\sigma$-invariant open pro-$p$-subgroup $\Gamma \subset\CG(\BQ_p)$, one can show that the subgroup generated by any two non-commuting odd elements of $\Gamma$ is an infinite strong Schur $\sigma$-group of Zassenhaus type $(3,3)$. 
\end{Rem}

Now we can establish the main result of this paper:

\begin{Thm}\label{I33SG4}
The set of $\sigma$-isomorphism classes of infinite weak Schur $\sigma$-groups of Zassenhaus type (3,3) is a closed subset of $\Sch$ of measure~$0$.
\end{Thm}

\begin{Proof}
For any weak Schur $\sigma$-group $G$, being of Zassenhaus type (3,3) depends only on the finite quotient $G/D_4(G)$. By the construction of the topology on~$\Sch$, the subset $\Sch_{(3,3)}$ of $\sigma$-isomorphism classes of weak Schur $\sigma$-groups of Zassenhaus type (3,3) is therefore open and closed in~$\Sch$. On the other hand, Proposition \ref{FinSchurOpenClosed} implies that the subset of $\sigma$-isomorphism classes of finite weak Schur $\sigma$-groups of Zassenhaus type (3,3) is an open subset. The set 
of $\sigma$-isomorphism classes of infinite weak Schur $\sigma$-groups of Zassenhaus type (3,3) is therefore closed in~$\Sch$.

Here the isomorphism classes of groups which are also strong form an at most countable set by Proposition \ref{InfStrSchSig33Fin}, which therefore has measure $0$ by Corollary \ref{CountInfNull}.
The remaining isomorphism classes form the subset $\Sch_{(3,3)}\setminus\Sch^\strong$, whose measure is~$0$ by Theorem~\ref{StrongThm}. Thus the union of these sets has measure $0$, and we are done.
\end{Proof}

\section{$p$-Tower groups}
\label{pTG}

Consider an imaginary quadratic field $K$ and let $K_p$ denote its maximal unramified pro-$p$-extension. Then $K_p$ is Galois over~$\BQ_p$ and its Galois group is the semidirect product of the pro-$p$-group $G_K := \Gal(K_p/K)$ with a group of order~$2$ generated by complex conjugation. Thus $G_K$ is a $\sigma$-pro-$p$-group. By Koch-Venkov \cite[\S1]{KochVenkov1974} and \cite[Prop.\,11.2]{PinkRubio2025}
this is a strong Schur $\sigma$-group. Following McLeman \cite{McLeman2008} we call it the \emph{$p$-tower group} associated to~$K$.

\medskip
In \cite[\S12]{PinkRubio2025} we explain that the probability measure $\mu_\infty$ on $\Sch$ is meant as a heuristic to describe statistical properties of $p$-tower groups. A particularly interesting question is when and how often $G_K$ is finite. Though much is known about that, the problem is not yet settled completely. 

\medskip
The present paper was motivated by the ``if'' part of McLeman's conjecture \cite[Conj.\,2.9]{McLeman2008}, which says that any $p$-tower group of Zassenhaus type (3,3) should be finite. 
Theorem \ref{I33SG4} yields some credence for this in the case $p>3$. Indeed, if our heuristic is correct, it implies that McLeman's statement holds on average, namely that the proportion of finite groups $G_K$ of Zassenhaus type (3,3) tends to~$1$. Of course, this does not rule out a sparse infinite set of exceptions.


\medskip
But recall that the Fontaine-Mazur Conjecture \cite[Conj.\,5b]{FontaineMazur1995} predicts that $G_K$ does not have an infinite $p$-adic analytic quotient. In particular this implies that $G_K$ itself cannot be an infinite $p$-adic analytic group (compare Hajir \cite[Conj.\,2]{Hajir1997}). But if $G_K$ is of Zassenhaus type (3,3) and $p>3$, Proposition \ref{InfStrong33ThenA1} does prove that $G_K$ is a $p$-adic analytic group. The Fontaine-Mazur conjecture therefore directly implies that $G_K$ is finite, yielding the ``if'' part of McLeman's conjecture in full in the case $p>3$.



\end{document}